\documentclass[11pt]{article}
\usepackage{amssymb, amsmath,fullpage}
\numberwithin{equation}{section}
\newtheorem {thm}{Theorem}[section]
\newtheorem {lem}[thm]{Lemma}

\newtheorem {prop}[thm]{Proposition}

\newtheorem {rem}[thm]{Remark}
\newcommand{\la}{\langle}
\newcommand{\ra}{\rangle}
\newcommand{\proof}{\par\noindent{\em Proof.\ }}

\newcommand{\qed}{\hfill $\square$\par\smallskip}
\title{\bf\boldmath 
Evaluation of the $BC_n$ elliptic Selberg integral 
\\
via the fundamental invariants 
}
\author{
{\sc Masahiko Ito}\thanks{
School of Science and Technology for Future Life,
Tokyo Denki University, Tokyo 120-8551, Japan
}\ \ 
and
{\sc Masatoshi Noumi}\thanks{
Department of Mathematics, Kobe University, 
Rokko, Kobe 657-8501, Japan
}
}
\date{}
\begin{document}
\maketitle
\begin{abstract}
We give an alternative proof of the evaluation formula for 
the elliptic Selberg integral of type $BC_n$
as an application of the fundamental $BC_n$-invariants. 
\end{abstract}


\section{Introduction}\label{sec:intro}
The evaluation formula of the $BC_n$ elliptic Selberg integral was proposed 
for the first time by van Diejen and Spiridonov \cite{vDS}. 
Namely, 
under the balancing condition $a_1\cdots a_6t^{2n-2}=pq$, 
\begin{equation}\label{eq:BCSum-0}
\begin{split}
&\frac{1}{(2\pi\sqrt{-1})^n}
\int_{\mathbb{T}^n}
\prod_{i=1}^{n}
\frac{
\prod_{m=1}^{6}\Gamma(a_mz_i^{\pm1};p,q)
}{
\Gamma(z_i^{\pm2};p,q)
}
\prod_{1\le j<k\le n}
\frac{
\Gamma(tz_j^{\pm1}z_k^{\pm1};p,q)
}{
\Gamma(z_j^{\pm1}z_k^{\pm1};p,q)
}
\frac{dz_1\cdots dz_n}{z_1\cdots z_n}
\\
&\qquad\qquad=
\frac{2^n n!}{(p;p)_\infty^n(q;q)_\infty^n}
\prod_{i=1}^{n}\bigg(
\frac{\Gamma(t^i;p,q)}{\Gamma(t;p,q)}
\prod_{1\le j<k\le 6}\Gamma(t^{i-1}a_ja_k;p,q)
\bigg),
\end{split}
\hspace{-6pt}
\end{equation}
where $a_1,\ldots,a_6$, $t$ are complex parameters 
with $|a_m|<1$ ($m=1,\ldots,6$), $|t|<1$, and 
$\mathbb{T}^n$ stands for the $n$-dimensional torus.  
(Here  $\Gamma(z;p,q)$ denotes the Ruijsenaars elliptic gamma function, 
and the double-signs indicate a product of all possible factors.) 
In the paper \cite{vDS}, they outlined a way of proof for \eqref{eq:BCSum-0} following Anderson's method \cite{An}, 
which is known as  a typical derivation for the evaluation formula of the Selberg integral \cite{Se} 
via the other multi-dimensional integral \cite{Di} called 
Dixon--Anderson integral in \cite{F,IF1}. 
The proof outlined in \cite{vDS} was eventually completed by Rains \cite{R}, 
proving the elliptic counterpart of 
the evaluation of the Dixon--Anderson integral
\begin{equation*}%
\begin{split}
&\frac{1}{(2\pi\sqrt{-1})^n}
\int_{\mathbb{T}^n}
\prod_{i=1}^{n}
\frac{
\prod_{m=1}^{2n+4}\Gamma(a_mz_i^{\pm1};p,q)
}{
\Gamma(z_i^{\pm2};p,q)
}
\prod_{1\le j<k\le n}
\frac{
1
}{
\Gamma(z_j^{\pm1}z_k^{\pm1};p,q)
}
\frac{dz_1\cdots dz_n}{z_1\cdots z_n}
\\
&\qquad\qquad=
\frac{2^n n!}{(p;p)_\infty^n(q;q)_\infty^n}
\prod_{1\le j<k\le 2n+4}\Gamma(a_ja_k;p,q),
\end{split}
\end{equation*}
whose alternative proof was given by Spiridonov \cite{S}.

Besides Anderson's method, several derivations are known for 
the evaluation formula of the Selberg integral.  
Aomoto \cite{Ao} gave an alternative proof 
by characterizing 
the integral as a solution of a difference equation with some specific boundary condition (see also \cite{IF2} for the $q$-integral case). 
The aim of this paper is to give an alternative proof for the $BC_n$ 
elliptic Selberg integral \eqref{eq:BCSum-0}, following Aomoto's method 
as is outlined below. 
Denoting by $I(a_1,\ldots,a_6)$ the left-hand side of \eqref{eq:BCSum-0}, 
we first prove that, 
under the balancing condition $a_1\cdots a_6t^{2n-2}=pq$, 
this integral satisfies the system of 
$q$-difference equations 
\begin{equation}\label{eq:qDE-0}
I(a_1,\ldots,a_5,a_6)
=I(a_1,\ldots,qa_k,\ldots,a_5,q^{-1}a_6)
\prod_{i=1}^{n}\prod_{1\le m\le 5\atop
m\ne k}
\frac{\theta(q^{-1}a_ma_6t^{i-1};p)}{\theta(a_ma_kt^{i-1};p)}
\hspace{-10pt}
\end{equation}
for $k=1,\ldots,5$.  
Setting 
\begin{equation*}
\widetilde{\Psi}(z)=
\prod_{i=1}^{n}
\frac{
\Gamma(pa_6z_i^{\pm1};p,q)
\prod_{m=1}^{5}\Gamma(a_mz_i^{\pm1};p,q)
}{
\Gamma(z_i^{\pm2};p,q)
}
\prod_{1\le j<k\le n}
\frac{
\Gamma(tz_j^{\pm1}z_k^{\pm1};p,q)
}{
\Gamma(z_j^{\pm1}z_k^{\pm1};p,q)
}, 
\end{equation*}
we use the notation 
\begin{equation*}
\la \varphi(z)\ra=\int_{\mathbb{T}^n}\varphi(z)\widetilde{\Psi}(z)\frac{dz_1\cdots dz_n}{z_1\cdots z_n}
\end{equation*}
for any meromorphic function $\varphi(z)$ on $(\mathbb{C}^\ast)^n$. 
Then 
the difference equation \eqref{eq:qDE-0} of the case $k=1$ is equivalent to the equality 
\begin{equation}
\label{eq:<E0n>}
\la E_n(a_1,a_6;z)\ra 
=
\la E_0(a_1,a_6;z)\ra 
\prod_{i=1}^{n}\left(
\frac{a_1^3\theta(a_6a_1^{-1}t^{i-1};p)}
{a_6^3\theta(a_1a_6^{-1}t^{i-1};p)}
\prod_{m=2}^{5}
\frac{\theta(a_ma_6t^{i-1};p)}{\theta(a_ma_1t^{i-1};p)}
\right)
\hspace{-10pt}
\end{equation}
under the balancing condition $a_1\cdots a_6t^{2n-2}=1$, 
where  
\begin{equation*}
E_0(a,b;z)=
\prod_{i=1}^{n}\frac{\theta(az_{i}^{\pm 1};p)}
{\theta(a(bt^{i-1})^{\pm 1};p)},
\quad
E_n(a,b;z)=
\prod_{i=1}^{n}\frac{\theta(bz_{i}^{\pm 1};p)}
{\theta(b(at^{i-1})^{\pm 1};p)}. 
\end{equation*}
The idea of Aomoto's method is 
to introduce appropriate intermediate functions 
which interpolate equation \eqref{eq:<E0n>}.  
We now define a set of holomorphic symmetric functions by 
\begin{equation}
\label{eq:E2-explicit02}
E_{r}(a,b;z)=
\hspace{-10pt}
\sum_{1\le i_1<\cdots<i_{r}\le n\atop 1\le j_1<\cdots<j_{n-r}\le n}
\prod_{k=1}^{r}\frac{\theta(b t^{i_k-k}z_{i_k}^{\pm 1};p)}{\theta(b t^{i_k-k}(at^{k-1})^{\pm 1};p)}
\prod_{l=1}^{n-r}\frac{\theta(at^{j_l-l}z_{j_l}^{\pm 1};p)}{\theta(at^{j_l-l}(bt^{l-1})^{\pm 1};p)}
\hspace{-10pt}
\end{equation}
for $r=0,1,\ldots,n$, 
where the summation is taken over all pairs of sequences 
$1\le i_1<\cdots<i_r\le n$ and $1\le j_1<\cdots<j_{n-r}\le n$ such that 
$\{i_1,\ldots,i_r\}\cup\{j_1,\ldots,j_{n-r}\}=\{1,2,\ldots,n\}$.
Under the condition $a_1\cdots a_6t^{2n-2}=1$, 
one can show that the following recurrence relations hold:
\begin{equation}
\label{eq:2-term}
\la E_{r}(a_1,a_6;z)\ra=C_{r}\la E_{r-1}(a_1,a_6;z)\ra
\quad (r=1,\ldots,n), 
\end{equation}
where the coefficients $C_r$ are given by 
\begin{equation*}\label{eq:ci}
C_{r}=
-\frac{
a_1^2t^{2r-2}\,
\theta(t^{n-r+1};p)
\theta(a_6a_1^{-1}t^{n-r+1};p)
\theta(a_1a_6^{-1}t^{2r-n};p)
}
{
a_6^2t^{2n-2r}\,
\theta(t^r;p)
\theta(a_6a_1^{-1}t^{n-2r+2};p)
\theta(a_1a_6^{-1}t^r;p)
}\!
\prod_{m=2}^5
\dfrac{\theta(a_ma_6t^{n-r};p)}{\theta(a_ma_1t^{r-1};p)}.
\end{equation*}
Using \eqref{eq:2-term} repeatedly, we immediately obtain \eqref{eq:<E0n>}. 
We call these $E_{r}(a,b;z)$ the {\em fundamental invariants} of type $BC_n$, which thus play an essential role in this paper. 
The fundamental invariants \eqref{eq:E2-explicit02} are given as 
a special case of the {\em Lagrange interpolation functions of type $BC_n$} 
in the context of the connection problem among the independent cycles for the $BC_n$ Jackson integral, see \cite[Example 2 of Theorem 1.4]{IN15}.
See also \cite{IN} for details of the fundamental invariants \eqref{eq:E2-explicit02}. 
We remark that our fundamental invariants $E_r(a,b;z)$ 
are essentially the interpolation theta functions of Coskun--Gustafson \cite{CG2007} and Rains \cite{Ra2006} 
attached to single columns of partitions. 
In fact, $E_r(a,b;z)$ are compared explicitly with the functions of \cite{CG2007} and \cite{Ra2006}, respectively, 
as explained in \cite[Introduction]{IN}. 
Also, the key equation \eqref{eq:2-term} is essentially the same as \cite[Theorem 4.1]{IN}
which we proved in the context of a $BC_n$ elliptic summation formula. 
It should be mentioned that van Diejen--Spiridonov \cite{vDS} already pointed out that 
the integral \eqref{eq:BCSum-0} implies the $BC_n$ elliptic summation formula via residue calculus.

Note that the integral \eqref{eq:BCSum-0} with $p=0$ is known as Gustafson's contour $q$-integral \cite{Gu}, 
which is the Nassrallah--Rahman integral in the case $n=1$ \cite{NR}. 
Aomoto's method using the fundamental invariants \eqref{eq:E2-explicit02} with $p=0$ 
leads us to the recurrence relations for the Gustafson's contour $q$-integral. 
This fact was previously discussed in \cite[Corollary 5.2 and Eq.\,(5.3)]{I}.

In order to establish the evaluation formula \eqref{eq:BCSum-0}, 
we need to investigate further 
the boundary condition for the difference 
equations \eqref{eq:qDE-0}; 
the precise arguments will be given later in Section \ref{section:4}. 

\par\medskip
This paper is organized as follows. After defining basic terminology in Section \ref{section:1}, 
we first discuss the system 
of $q$-difference equations \eqref{eq:qDE-0}
in Section \ref{section:2}. 
In Section \ref{section:3} we study the analytic continuation of the integral \eqref{eq:BCSum-0} 
as a meromorphic function of the parameters $a_1,\ldots,a_5$ in a specific domain. 
We use this argument to show that the integral \eqref{eq:BCSum-0} is expressed as a product of elliptic gamma functions up to a constant. 
Section \ref{section:4} is devoted to obtaining the boundary condition for \eqref{eq:qDE-0}
through asymptotic analysis of 
the contour integral \eqref{eq:BCSum-0} 
as $a_2\to a_1^{-1}$ (i.e. $a_1a_2\to 1$). 
This condition determines the explicit value of the constant, which was indefinite at the time of Section \ref{section:3}. 
In the case of elliptic hypergeometric integrals, 
we often meet some strict restraints on parameters, 
which do not permit us to consider the asymptotic behavior like $a_i \to 0$ or $\infty$ as we usually do in the rational or trigonometric ($q$-analog) cases. 
Thus our treatment of the boundary condition might look totally different from that of the $q$-analog case. 
It should  be noted, however, that 
our method to analyze such a situation as $a_1a_2\to 1$
is also applicable to the case $p=0$ of the integral \eqref{eq:BCSum-0}, 
thus providing a novel insight even for the evaluation of contour $q$-integrals. 

\section{$BC_n$ elliptic Selberg integral}\label{section:1}
Throughout this paper we denote by 
$\Gamma(u;p,q)$ $(u\in\mathbb{C}^\ast)$ 
the Ruijsenaars elliptic gamma function defined by 
\begin{equation*}
\Gamma(u;p,q)=\frac{(pqu^{-1};p,q)_\infty}{(u;p,q)_\infty},\quad
(u;p,q)_\infty=\prod_{\mu,\nu=0}^{\infty}(1-p^\mu q^\nu u)\quad(|p|<1,|q|<1). 
\end{equation*}
Note that $\Gamma(u;p,q)$ satisfies
\begin{equation}
\label{eq:RGamma}
\Gamma(qu;p,q)=\theta(u;p)\Gamma(u;p,q),\quad
\Gamma(pq/u;p,q) =\frac{1}{\Gamma(u;p,q)}. 
\end{equation}
We consider the meromorphic function 
\begin{equation*}
\Psi(z)
=\prod_{i=1}^{n} 
\frac{\prod_{m=1}^{6}\Gamma(a_mz_i^{\pm 1};p,q)}
{\Gamma(z_i^{\pm2};p,q)}
\prod_{1\le j<k\le n}
\frac{\Gamma(tz_j^{\pm1}z_k^{\pm1};p,q)}
{\Gamma(z_j^{\pm1}z_k^{\pm1};p,q)} 
\end{equation*}
in $z=(z_1,\ldots,z_n)\in(\mathbb{C}^\ast)^n$ with complex parameters $a_1,\ldots,a_6,t\in\mathbb{C}^\ast$, 
assuming 
throughout that $|t|<1$.   
We also use the notation $\Psi(a_1,\ldots,a_6;z)$ for $\Psi(z)$ 
when we need to make 
the dependence on the parameters $a_1,\ldots,a_6$ explicit. 
For this function $\Psi(z)$, 
we investigate the multiple integral 
\begin{equation*}
I=
\int_{\sigma}\Psi(z)\varpi(z),\quad 
\varpi(z)=\frac{1}{(2\pi\sqrt{-1})^n}\frac{dz_1\cdots dz_n}{z_1\cdots z_n}
\end{equation*}
over an $n$-cycle $\sigma$.  
Since $\Psi(z)$ is expressed as 
\begin{equation*}
\begin{split}
\Psi(z)
&=
\prod_{i=1}^{n}\,
(1-z_i^{\pm2})(pz_i^{\pm 2};p)_\infty(qz_i^{\pm 2};q)_\infty
\\
&\quad\times
\prod_{1\le j<k\le n}
(1-z_j^{\pm1}z_k^{\pm1})(pz_j^{\pm 1}z_k^{\pm 1};p)_\infty(qz_j^{\pm 1}z_k^{\pm 1};q)_\infty
\\
&\quad\times
\prod_{i=1}^{n}\prod_{m=1}^{6}
\frac{(pqa_m^{-1}z_i^{\pm1};p,q)_\infty}{(a_mz_i^{\pm1};p,q)_\infty}
\prod_{1\le j<k\le n}
\frac{(pqt^{-1}z_j^{\pm1}z_k^{\pm1};p,q)_\infty}{(tz_j^{\pm1}z_k^{\pm1};p,q)_\infty}, 
\end{split}
\end{equation*}
we see that $\Psi(z)$ has poles possibly along the divisors
\begin{equation*}
\begin{split}
z_i^{\pm 1}&=a_m p^\mu q^\nu \quad
\ (1\le i\le n;\ m=1,\ldots,6;\ \mu,\nu=0,1,2,\ldots)
\\
z_j^{\pm1}z_k^{\pm1}&=t\,p^\mu q^{\nu}\qquad(1\le j<k\le n;\  \mu,\nu=0,1,2,\ldots). 
\end{split}
\end{equation*}
Also, 
regarded as a function of $z_i$ ($i=1,\ldots,n$), $\Psi(z)$ 
has poles possibly at 
\begin{equation*}
p^{\mu}q^{\nu}a_m,\quad 
p^{-\mu}q^{-\nu}a_m^{-1},\quad 
p^{\mu}q^{\nu}t z_j^{\pm1},\quad
p^{-\mu}q^{-\nu}t^{-1}z_j^{\pm1},
\end{equation*}
where $1\le m\le 6$, $1\le j\le n$, $j\ne i$ and 
$\mu,\nu=0,1,2,\ldots$.  
If the parameters satisfy the condition $|a_1|<1,\ldots,|a_6|<1$, 
then $\Psi(z)$ is holomorphic 
in a neighborhood of the $n$-dimensional torus
\begin{equation*}
\mathbb{T}^n=\big\{\ z=(z_1,\ldots,z_n)\in(\mathbb{C}^\ast)^n\ \big\vert\ 
|z_i|=1\quad(i=1,\ldots,n) \ \big\}, 
\end{equation*}
and hence the integral 
\begin{equation*}
I(a_1,\ldots,a_6)=\int_{\mathbb{T}^n}\Psi(a_1,\ldots,a_6;z)\varpi(z)
\end{equation*}
defines a holomorphic function on the domain 
\begin{equation}
\label{eq:defU}
U=\big\{\  (a_1,\ldots,a_6)\in(\mathbb{C}^\ast)^6\ \big\vert\ 
|a_m|<1\quad(m=1,\ldots,6) \ \big\}\subset(\mathbb{C}^\ast)^6.   
\end{equation}  
This function can be continued to a 
holomorphic function on a larger domain by 
replacing $\mathbb{T}^n$ with an appropriate 
$n$-cycle depending on the parameters 
$(a_1,\ldots,a_6)$.  
We give below a remark on analytic continuation of this sort. 

\par\medskip
For each $(a_1,\ldots,a_6)\in(\mathbb{C}^\ast)^6$, 
we define two subsets $S_0$, $S_\infty$ 
of $\mathbb{C}^\ast$ by
\begin{equation*}
\begin{split}
S_0&=\big\{\ p^{\mu}q^{\nu}a_m\ \big\vert\  
1\le m\le 6;\ \mu,\nu\in\mathbb{N}\ \big\},
\\ 
S_\infty&=\big\{\ p^{-\mu}q^{-\nu}a_m^{-1}\ \big\vert\  
1\le m\le 6;\ \mu,\nu\in\mathbb{N}\ \big\}, 
\end{split}
\end{equation*}
where $\mathbb{N}=\{0,1,2,\ldots\}$, 
and suppose that $S_0\cap S_\infty=\phi$.  
Assuming that $|t|<r^2$ for some $r\in(0,1)$, 
we choose a circle 
\begin{equation*}
C_\rho(0)=\big\{\ u\in\mathbb{C}^\ast\ \big\vert\ 
|u|=\rho\ \big\},\quad \rho\in[r,r^{-1}], 
\end{equation*}
which does not intersect with $S_0\cup S_\infty$.  
Then we define a cycle $C$ in $\mathbb{C}^\ast$
by
\begin{equation*}
C=C_\rho(0)+
\sum_{c\in S_0; |c|>\rho}\,C_\varepsilon(c)-
\sum_{c\in S_\infty; |c|<\rho}\,C_\varepsilon(c),
\end{equation*}
where $C_\varepsilon(c)$ denotes a sufficiently small circle 
around $c$.  
Note that, if $|a_m|<1$ ($m=1,\ldots,6$), then $C$ 
is homologous to the unit circle. 
We now assume that  $|a_m|<r^{-1}$ ($m=1,\ldots,6$).   
Then such a cycle $C$ can be taken inside the annulus 
$A_r=\{ u\in\mathbb{C}^\ast\ \vert\ r\le |u|\le r^{-1}\}$.  
Since $|t|<r^2$, the meromorphic function 
$\Psi(z)$ is holomorphic in an neighborhood of the 
$n$-cycle $C^n=C\times\cdots\times C$.  
Hence, the integral 
\begin{equation*}
I=\int_{C^n}\Psi(z)\varpi(z)
\end{equation*}
is well defined, and does not depend on the choice 
of $\rho\in [r,r^{-1}]$.  This implies the following 
lemma on analytic continuation.

\begin{lem}\label{lem:analcont}
Suppose that $|t|<r^2$ for some real number
 $r\in(0,1]$.  Then the holomorphic function $I(a_1,\ldots,a_6)$ 
on the domain $U$ of \eqref{eq:defU}
can be continued to a holomorphic function on 
\begin{equation}\label{eq:defUext}
\bigg\{(a_1,\ldots,a_6)\in(\mathbb{C}^\ast)^6
\ \bigg\vert\ 
{
|a_m|<r^{-1}\ (1\le m\le 6),\ 
\atop
a_ka_l\notin p^{-\mathbb{N}}q^{-\mathbb{N}}
\ (1\le k,l\le 6)
}
\ \bigg\}.
\end{equation}
\qed
\end{lem}

As can be seen in \eqref{eq:BCSum-0}, under the balancing condition 
$a_1\cdots a_6t^{2n-2}=pq$, this function 
$I(a_1,\ldots,a_6)$ is eventually continued to a meromorphic function on 
a hypersurface in $(\mathbb{C}^\ast)^6$ with poles along 
the divisors 
\begin{equation*}
p^{\mu}q^{\nu}t^{i-1}a_ka_l=1\quad(k,l\in\{1,\ldots,6\}; i=1,\ldots,n; 
\mu,\nu=0,1,2,\ldots). 
\end{equation*}

\section{$q$-Difference equations with respect to the parameters}\label{section:2}
In this section we derive a system of 
$q$-difference equations for 
the integral $I(a_1,\ldots,a_6)$ on the basis of the arguments in 
\cite{IN}.  
Our goal is to establish the following proposition.  
\begin{prop}
\label{prop:qDE} 
Suppose that $|p|<|t|^{2n-2}$.  
Under the balancing condition $a_1\cdots a_6t^{2n-2}=pq$, 
the integral $I(a_1,\ldots,a_6)$ satisfies the system of 
$q$-difference equations 
\begin{equation}
\label{eq:qDE}
I(a_1,\ldots,a_5,a_6)=I(a_1,\ldots,qa_k,\ldots,a_5,q^{-1}a_6)
\prod_{i=1}^{n}\prod_{1\le m\le 5
\atop m\ne k}
\!\!
\frac{\theta(q^{-1}a_ma_6t^{i-1};p)}{\theta(a_ma_kt^{i-1};p)}
\hspace{-5pt}
\end{equation}
for $k=1,\ldots,5$, 
provided that $|a_1|<1,\ldots,|a_5|<1$ and $|a_6|<|q|$.  
\end{prop}
Note that the condition $|a_6|<|q|$ is equivalent 
to $|a_1\cdots a_5|>|p|/|t|^{2n-2}$ under the balancing 
condition.  
We need to assume that $p$ is sufficiently 
small as specified above to guarantee 
that \eqref{eq:qDE} 
holds in a nonempty region. 

\par\medskip
In order to make use of the arguments of \cite{IN}, we modify 
$\Psi(z)$ as
\begin{equation*}
\begin{split}
\widetilde{\Psi}(z)&=\Psi(a_1,\ldots,a_5,pa_6;z)
\\
&=
\prod_{i=1}^{n}
\frac{
\Gamma(pa_6z_i^{\pm1};p,q)
\prod_{m=1}^{5}\Gamma(a_mz_i^{\pm1};p,q)
}{
\Gamma(z_i^{\pm2};p,q)
}
\prod_{1\le j<k\le n}
\frac{
\Gamma(tz_j^{\pm1}z_k^{\pm1};p,q)
}{
\Gamma(z_j^{\pm1}z_k^{\pm1};p,q)
}
\\
&=
\prod_{i=1}^{n}
\frac{
\prod_{m=1}^{5}\Gamma(a_mz_i^{\pm1};p,q)
}{
\Gamma(qa_6^{-1}z_i^{\pm1};p,q)
\Gamma(z_i^{\pm2};p,q)
}
\prod_{1\le j<k\le n}
\frac{
\Gamma(tz_j^{\pm1}z_k^{\pm1};p,q)
}{
\Gamma(z_j^{\pm1}z_k^{\pm1};p,q)
}. 
\end{split}
\end{equation*} 
This function 
$\widetilde{\Psi}(z)$ coincides with the meromorphic function 
$\widetilde{\Phi}(z)$ in \cite{IN} up to multiplication 
by a $q$-periodic function in all variables $z_1,\ldots,z_n$. 
Namely one has 
\begin{equation*}
\begin{split}
\frac{T_{q,z_i}\widetilde{\Psi}(z)}{\widetilde{\Psi}(z)}
&=-
\frac{(q^{-1}z_i^{-1})^2\theta(q^{-2}z_i^{-2};p)}{z_i^2\,\theta(z_i^2;p)}
\prod_{m=1}^{6}\frac{\theta(a_mz_i;p)}
{\theta(q^{-1}a_mz_i^{-1};p)}
\\
&\quad\times\,
\prod_{1\le k\le n
\atop k\ne i}
\frac{\theta(tz_iz_k^{\pm1};p)
\theta(q^{-1}z_i^{-1}z_k^{\pm1};p)}
{\theta(q^{-1}tz_i^{-1}z_k^{\pm 1};p)
\theta(z_iz_k^{\pm1};p)} 
\end{split}
\end{equation*}
for $i=1,\ldots,n$, where $T_{q,z_i}$ stands for the $q$-shift operator 
in $z_i$:  
\begin{equation*}
T_{q,z_i}f(z_1,\ldots,z_n)=f(z_1,\ldots,qz_i,\ldots,z_n). 
\end{equation*}
As to the parameters $a_1,\ldots,a_6$, 
one has 
\begin{equation}
\label{eq:qshifts}
\begin{split}
\frac{T_{q,a_m}\widetilde{\Psi}(z)}{\widetilde{\Psi}(z)}
&=\prod_{i=1}^{n}\theta(a_mz_i^{\pm1};p)\quad(1\le m\le 5),
\\
\frac{T_{q,a_6}\widetilde{\Psi}(z)}{\widetilde{\Psi}(z)}
&=a_6^{-2n}\prod_{i=1}^{n}\theta(a_6z_i^{\pm1};p).  
\end{split}
\end{equation}
In this paper we use the notation of expectation values 
to refer to the integral
\begin{equation*}
\la \varphi(z)\ra=\int_{\mathbb{T}^n}\varphi(z)\widetilde{\Psi}(z)\varpi(z)
\end{equation*}
for any meromorphic function $\varphi(z)$ on $(\mathbb{C}^\ast)^n$ 
such that 
$\varphi(z)\widetilde{\Psi}(z)$ is holomorphic in a neighborhood of the 
$n$-dimensional torus $\mathbb{T}^n$.  
If we set 
\begin{equation*}
\nabla_{q,z_i}\varphi(z)=\varphi(z)-\dfrac{T_{q,z_i}\widetilde{\Psi}(z)}{\widetilde{\Psi}(z)}
T_{q,z_i}\varphi(z)\qquad(i=1,\dots,n)
\end{equation*}
as in \cite{IN}, one has 
\begin{equation*}
\la \nabla_{q,z_i}\varphi(z)\ra=0\qquad(i=1,\ldots,n)
\end{equation*}
for any meromorphic function $\varphi(z)$ such that 
$\varphi(z)\widetilde{\Psi}(z)$ is holomorphic in a neighborhood 
of the compact set
\begin{equation}\label{eq:Tqi}
|q|\le |z_i|\le 1, \quad |z_j|=1\quad (1\le j\le n;\ j\ne i).
\end{equation} 
In fact, by the Cauchy 
theorem one has
\begin{equation*}
\int_{\mathbb{T}^n}\widetilde{\Psi}(z)\varphi(z)\varpi(z)=
\int_{\mathbb{T}^n}T_{q,z_i}(\widetilde{\Psi}(z)\varphi(z))\varpi(z)
\quad(i=1,\ldots,n). 
\end{equation*}

We set 
\begin{equation*}
K(a_1,\ldots,a_5,a_6)=I(a_1,\ldots,a_5,pa_6)=\la 1\ra 
\end{equation*}
assuming that 
$|a_m|<1$ $(m=1,\ldots,5)$, $|pa_6|<1$. 
Then from \eqref{eq:qshifts} we have 
\begin{equation}\label{eq:KK}
\begin{split}
K(qa_1,a_2,\ldots,a_6)&
=\la \prod_{i=1}^{n}\theta(a_1z_i^{\pm1};p)\ra 
=\la E_0(a_1,a_6;z)\ra \prod_{i=1}^{n}\theta(a_1(a_6t^{i-1})^{\pm1};p),
\\
K(a_1,\ldots,a_5,qa_6)
&
=a_6^{-2n}\la \prod_{i=1}^{n}\theta(a_6z_i^{\pm1};p)\ra 
=\la E_n(a_1,a_6;z)\ra \prod_{i=1}^{n}a_6^{-2}\theta(a_6(a_1t^{i-1})^{\pm1};p), 
\end{split}
\end{equation}
where $E_{r}(a_1,a_6;z)=E_{r}^{(n)}(a_1,a_6;z)$ ($r=0,1,\ldots,n$) denote 
the fundamental $BC_n$-invariants \eqref{eq:E2-explicit02}; 
for the basic properties of these functions, we refer the reader to 
\cite[Section 3]{IN}. 
On the other hand, 
as for the function 
\begin{equation*}
\varphi_{r,i}(z)=F_i^{-}(z)E^{(n-1)}_{r-1}(a_1,a_6;z_1,\ldots,z_{i-1},z_{i+1},\ldots,z_{n})
\hspace{10pt} (1\le i\le n; 1\le r\le n)
\end{equation*}
of \cite[Section 4]{IN}, 
where
\begin{equation*}
F_{i}^-(z)=
\frac{\prod_{m=1}^{6}
\theta(a_mz_i^{-1};p)}{z_i^{-2}\, \theta(z_i^{-2};p)}
\prod_{1\le j\le n
\atop j\ne i}
\frac{\theta(tz_i^{-1}z_j^{\pm 1};p)}{\theta(z_i^{-1}z_j^{\pm 1};p)},
\end{equation*}
one can verify that 
$\varphi_{r,i}(z)\widetilde{\Psi}(z)$ is holomorphic in a neighborhood of \eqref{eq:Tqi}.  
In fact, in the product $F_i^{-}(z)\widetilde{\Psi}(z)$, 
all possible poles of each of the two functions $F_i^{-}(z)$, 
$\widetilde{\Psi}(z)$ 
\begin{equation*}
\begin{split}
&p^{\mu} z_i^{-2}=1, \quad p^{\mu}z_iz_j^{\pm1}, 
\quad
p^{\mu}a_m z_i^{-1},\quad
p^{\mu} tz_i^{-1}z_j^{\pm1}=1
\\
&\quad(1\le j\le n;\ j\ne i;\ m=1,\ldots,6; \ \mu=0,1,2,\ldots)
\end{split}
\end{equation*}
relevant to this region  are eliminated 
by zeros of the other.  
Hence we have
\begin{equation*}
\la\nabla_{q,z_i}\varphi_{r,i}(z)\ra=0\qquad(i=1,\ldots,n).  
\end{equation*}
In the same way as we discussed in \cite[Theorem 4.1]{IN}, 
this formula implies the recurrence relation \eqref{eq:2-term}, 
and hence
\begin{equation*}
\la E_n(a_1,a_6;z)\ra =
\la E_0(a_1,a_6;z)\ra 
\prod_{i=1}^{n}\left(
\frac{a_1^3\theta(a_6a_1^{-1}t^{i-1};p)}
{a_6^3\theta(a_1a_6^{-1}t^{i-1};p)}
\prod_{m=2}^{5}
\frac{\theta(a_ma_6t^{i-1};p)}{\theta(a_ma_1t^{i-1};p)}
\right)
\end{equation*}
under the balancing condition $a_1\cdots a_6t^{2n-2}=1$.  
Combining this with \eqref{eq:KK} we obtain 
\begin{equation*}
\begin{split}
K(a_1,\ldots,a_5,qa_6)
&=
K(qa_1,a_2,\ldots,a_6)
\dfrac{a_1^{n}}{a_6^{3n}}
\prod_{i=1}^{n}
\prod_{m=2}^{5}
\frac{\theta(a_ma_6t^{i-1};p)}{\theta(a_ma_1t^{i-1};p)}
\\
&=
K(qa_1,a_2,\ldots,a_6)
(a_1\cdots a_6t^{2n-2})^n
\prod_{i=1}^{n}
\prod_{m=2}^{5}
\frac{\theta(pa_ma_6t^{i-1};p)}{\theta(a_ma_1t^{i-1};p)}
\\
&=
K(qa_1,a_2,\ldots,a_6)
\prod_{i=1}^{n}
\prod_{m=2}^{5}
\frac{\theta(pa_ma_6t^{i-1};p)}{\theta(a_ma_1t^{i-1};p)}.  
\end{split}
\end{equation*}
In terms of the function $I(a_1,\ldots,a_6)$, we conclude that
\begin{equation}\label{eq:qDE0}
I(a_1,\ldots,a_5,pqa_6)=
I(qa_1,a_2,\ldots,a_5,pa_6)
\prod_{i=1}^{n}
\prod_{m=2}^{5}
\frac{\theta(pa_ma_6t^{i-1};p)}{\theta(a_ma_1t^{i-1};p)} 
\end{equation}
under the conditions 
$a_1\cdots a_6 t^{2n-2}=1$, 
$|a_m|<1$ ($m=1,\ldots,5$) and $|pa_6|<1$.  
Hence, 
replacing $pa_6$ by $a_6$ in \eqref{eq:qDE0}
and changing the balancing condition accordingly, 
we have 
\begin{lem}\label{lem:qDE}
Under the conditions 
$a_1\cdots a_6t^{2n-2}=p$ 
and $|a_m|<1$ $(m=1,\ldots,6)$,  
one has 
\begin{equation}\label{eq:qDEA}
I(a_1,\ldots,a_5,qa_6)
=
I(a_1,\ldots,qa_k,\ldots,a_6)
\prod_{i=1}^{n}
\prod_{1\le m\le 5
\atop m\ne k}
\frac{\theta(a_ma_6t^{i-1};p)}{\theta(a_ma_kt^{i-1};p)} 
\end{equation}
for $k=1,\ldots,5$. 
\qed
\end{lem}
Further, 
replacing $a_6$ by $q^{-1}a_6$ we obtain Proposition 
\ref{prop:qDE}. 

\par\medskip
We now suppose that $a_1\cdots a_6t^{2n-2}=pq$, 
and regard $a_6=pq/a_1\cdots a_5t^{2n-2}$ 
as a function of $(a_1,\ldots,a_5)$.  
Then the integral $I(a_1,\ldots,a_6)$, regarded as 
a function of $(a_1,\ldots,a_5)$, 
is defined on the open subset 
\begin{equation*}
U_0=\big\{\ 
(a_1,\ldots,a_5)\in(\mathbb{C}^\ast)^5
\ \big\vert\ 
|a_1|<1,\ldots,|a_5|<1, |a_1\cdots a_5|>\frac{|p||q|}{|t|^{2n-2}}  
\ \big\}
\end{equation*}
of $(\mathbb{C}^\ast)^5$; we need to assume   
$|p||q|<|t|^{2n-2}$ in order to ensure that $U_0$ 
is not empty.  
We denote by 
\begin{equation*}
V_0=\big\{\ 
(a_1,\ldots,a_5)\in(\mathbb{C}^\ast)^5
\ \big\vert\ 
|a_1|<1,\ldots,|a_5|<1, |a_1\cdots a_5|>\frac{|p|}{|t|^{2n-2}} 
\ \big\} 
\end{equation*}
the nonempty open subset of $U_0$ where 
$I(a_1,\ldots,a_6)$ satisfies the $q$-difference equations
\eqref{eq:qDE}, assuming that $|p|<|t|^{2n-2}$. 

\section{Analytic continuation} \label{section:3}

The integral $I(a_1,\ldots,a_6)$, regarded as a holomorphic 
function in $(a_1,\ldots,a_5)\in U_0$, 
can be continued to a meromorphic function on 
$(\mathbb{C}^\ast)^5$.  
We prove this fact by means the $q$-difference 
equations \eqref{eq:qDE}.  
\par\medskip
In view of Proposition \ref{prop:qDE} 
we consider the meromorphic function
\begin{equation}\label{eq:defJ}
J(a_1,\ldots,a_6)=
\prod_{i=1}^{n}\,
\prod_{1\le j<k\le 6}\Gamma(a_ja_kt^{i-1};p,q).  
\end{equation}
Then it turns out that $J(a_1,\ldots,a_6)$ 
satisfies the same $q$-difference 
equations as \eqref{eq:qDE}.  
In fact, from \eqref{eq:RGamma} one has 
\begin{equation*}
J(a_1,\ldots,a_6)
=
J(a_1,\ldots,qa_k,\ldots,q^{-1}a_6)
\prod_{i=1}^{n}
\prod_{1\le m\le 5
\atop m\ne k}
\frac{\theta(q^{-1}a_ma_6;p)}
{\theta(a_ma_k;p)}. 
\end{equation*}
for $k=1,\ldots,5$.  
In the following we regard $J(a_1,\ldots,a_6)$ as 
a meromorphic function in $(a_1,\ldots,a_5)$ through 
$a_6=pq/a_1\cdots a_5t^{2n-2}$ as before.  

Noting that the integral 
$I(a_1,\ldots,a_6)$ is a holomorphic function on 
$U_0$, we consider the meromorphic function 
\begin{equation*}
\begin{split}
f(a_1,\ldots,a_6)
&=\frac{I(a_1,\ldots,a_6)}{J(a_1,\ldots,a_6)}
=\frac{I(a_1,\ldots,a_6)}
{\prod_{i=1}^{n}
\prod_{1\le j<k\le 6}\Gamma(a_ja_kt^{i-1};p,q)}
\\
&=
I(a_1,\ldots,a_6)\prod_{i=1}^{n}
\prod_{1\le j<k\le 6}
\frac
{(a_ja_kt^{i-1};p,q)_\infty}
{(pq/a_ja_kt^{i-1};p,q)_\infty}
\end{split}
\end{equation*}
on $U_0$.  
This ratio $f(a_1,\ldots,a_6)$ 
has poles possibly along the divisors
\begin{equation*}
t^{i-1}a_ja_k=p^{\mu+1} q^{\nu+1}\quad 
(i=1,\ldots,n;\ 1\le j<k\le 6;\ \mu,\nu=0,1,2,\ldots)  
\end{equation*}
in $U_0$. 
Also, 
$f(a_1,\ldots,a_6)$ 
is $q$-periodic with respect to $(a_1,\ldots,a_5)\in V_0$ 
in the sense that 
\begin{equation*}
f(a_1,\ldots,a_6)=f(a_1,\ldots,qa_k,\ldots,q^{-1}a_6)
\end{equation*}
for $k=1,\ldots,5$.  

\begin{lem}\label{lem:W0}
Suppose that $|p|<|q|^{\frac{25}{4}}|t|^{2n-2}$.
Then there exists an open subset 
$W_0\subset (\mathbb{C}^\ast)^5$
of the form
\begin{equation}\label{eq:W0def}
\begin{split}
W_0&=\big\{\ (a_1,\ldots,a_5)\in(\mathbb{C}^\ast)^5\ 
\big\vert\ 
sr<|a_m|<r\ (1\le m\le 5)\ \big\}
\hspace{15pt}
(0<r\le 1;\ 0<s<|q|)
\end{split}
\end{equation}
such that $W_0\subset V_0$ and that 
$f(a_1,\ldots,a_6)$ is holomorphic on $W_0$. 
\end{lem}
%
\proof
Under the assumption $|p|<|q|^{\frac{25}{4}}|t|^{2n-2}$, 
one can choose positive numbers $r,s$ such that 
\begin{equation*}
0<r<|q|^{\frac{1}{4}},\quad  
r^4|t|^{n-1}\le s<|q|,\quad |p|\le s^5r^5|t|^{2n-2}.  
\end{equation*}
Suppose that $(a_1,\ldots,a_5)\in W_0$.  Then 
$|a_1\cdots a_5|>s^5r^5\ge |p|/|t|^{2n-2}$ and hence 
$|a_6|<|q|$.  This means that $W_0\subset V_0$.  
Note also
$|a_6|=|pq/a_1\cdots a_5t^{2n-2}|>|p||q|/r^5|t|^{2n-2}$, 
and hence
\begin{equation*}
|p||q|/r^5|t|^{2n-2}<|a_6|<|q|.  
\end{equation*}
To show that $f(a_1,\ldots,a_6)$ is holomorphic in $W_0$ 
we verify
\begin{equation*}
|t^{n-1}a_ja_k|>|p||q|\quad(1\le j<k\le 6).  
\end{equation*}
In fact we have for $j=1,\ldots,5$, 
\begin{equation*}
|t^{n-1}a_ja_6|>|t|^{n-1}sr|p||q|/r^5|t|^{2n-2}=
|p||q|s/r^4|t|^{n-1}\ge |p||q|, 
\end{equation*}
and for $1\le j<k\le 5$, \\[8pt]
\hspace{115pt}$
|t^{n-1}a_ja_k|>s^2r^2|t|^{n-1}>s^5r^5|t|^{2n-2}\ge|p|> |p||q|. 
$\hfill $\square$

\begin{rem}\ \rm
If $|p|\le|q|^{10}|t|^{2n-2}$, 
one can simply take 
$r=|q|^{\frac{1}{2}}$ and 
$s=|q|^{\frac{3}{2}}$ 
for $f(a_1,\ldots,a_6)$ to be holomorphic on $W_0\subset V_0$. 
\end{rem}

\begin{thm}\label{thm:I=cJ}
Suppose that $|p|<|q|^{\frac{25}{4}}|t|^{2n-2}$. 
Under the condition $a_1\cdots a_6t^{2n-2}=pq$, 
the integral 
$I(a_1,\ldots,a_6)$, regarded as a holomorphic 
function in $(a_1,\ldots,a_5)\in U_0$, is continued to 
a meromorphic function on $(\mathbb{C}^\ast)^5$.   
Furthermore, it is expressed as
\begin{equation*}
I(a_1,\ldots,a_6)=c_n \prod_{i=1}^{n} 
\prod_{1\le j<k\le 6}\Gamma(a_ja_kt^{i-1};p,q)
\end{equation*}
for some constant $c_n\in\mathbb{C}$  
independent of $a_1,\ldots,a_6$.  
\end{thm}
\proof
By Lemma \ref{lem:W0}, there exists an open subset 
$W_0\subset(\mathbb{C}^\ast)^5$ of the form 
\eqref{eq:W0def} where $f(a_1,\ldots,a_6)$ 
is holomorphic and satisfies 
the $q$-difference equations 
\begin{equation}\label{eq:qDEf}
f(a_1,\ldots,a_6)=f(a_1,\ldots,q a_k,\ldots,q^{-1}a_6)
\qquad(k=1,\ldots,5). 
\end{equation}
for $(a_1,\ldots,a_5)\in W_0$.  
Note that $W_0$ 
is the product of $5$ copies of an annulus 
in which the ratio of the two radii is given by 
$s<|q|$.  Hence, by the $q$-difference equations 
\eqref{eq:qDEf}, the holomorphic function 
$f(a_1,\ldots,a_6)$ on $W_0$ is continued to a 
{\em holomorphic} function on the whole 
$(\mathbb{C}^\ast)^5$.  
It must be a constant, however, since the continued 
function $f(a_1,\ldots,a_6)$ is $q$-periodic with 
respect to the variables $a_1,\ldots,a_5$.  
If we denote this constant by $c_n$, we 
have $I(a_1,\ldots,a_6)=c_n J(a_1,\ldots,a_6)$ 
as a meromorphic function on $(\mathbb{C}^\ast)^5$.
\qed 

We compute the constant $c_n$ in the next section
by induction on the dimension $n$.  
Once this constant has been determined, we see 
that the statement above is valid for $|p|<1$ 
without any particular restriction. 

\section{Computation of the constant $c_n$ }\label{section:4}
In order to make the dimension explicit, 
we use below the notation 
$\Psi_n(z)$, $I_n(a_1,\ldots,a_6)$, $J_n(a_1,\ldots,a_6)$
for 
$\Psi(z)$, $I(a_1,\ldots,a_6)$, $J(a_1,\ldots,a_6)$
of the previous sections.  
As before, we assume that the parameters satisfy the 
balancing condition $a_1\cdots a_6t^{2n-2}=pq$, 
and regard $a_{6}=qp/a_1\cdots a_5t^{2n-2}$ as 
a function of $(a_1,\ldots,a_5)$.
By Theorem \ref{thm:I=cJ} 
we already know that two meromorphic functions 
$I_n(a_1,\ldots,a_6)$, $J_n(a_1,\ldots,a_6)$ are 
related by the formula 
\begin{equation}\label{eq:I=cJ}
I_n(a_1,\ldots,a_6)=c_n J_n(a_1,\ldots,a_6)
\end{equation}
provided that $|p|$ is sufficiently small.  
To determine the constant $c_n$, we investigate the 
behavior of the these two functions along the divisor $a_1a_2=1$. 

\par\medskip
We first consider the limit of $J_n(a_1,\ldots,a_6)$ as $a_2\to a_1^{-1}$. 
Noting that 
\begin{equation*}
\lim_{a_2\to a_1^{-1}}
(1-a_1a_2)
\Gamma(a_1a_2;p,q)=
\frac{1}{(p;p)_\infty(q;q)_\infty}, 
\end{equation*}
from \eqref{eq:defJ}
we have
\begin{equation*}
\begin{split}
&
\lim_{a_2\to a_1^{-1}}
(1-a_1a_2)J_n(a_1,\ldots,a_6)
\\&
=
\frac{\prod_{i=1}^{n-1}\Gamma(t^{i};p,q)}{(p;p)_\infty(q;q)_\infty}
\prod_{i=1}^{n}
\Bigg(
\prod_{m=3}^{6}
\Gamma(a_1^{\pm1}a_kt^{i-1};p,q)
\prod_{3\le j<k\le 6}{}
\Gamma(a_ja_kt^{i-1};p,q)
\Bigg),
\end{split}
\end{equation*}
where $a_6$ in the right-hand side should be understood as 
$a_6=pq/a_3a_4a_5t^{2n-2}$. 
Since 
$a_3a_4a_5a_6t^{2n-2}=pq$ in the limit, for any permutation 
$(i,j,k,l)$ of 
$(3,4,5,6)$ we have 
$(a_ia_jt^{n-1})(a_ka_lt^{n-1})=pq$, and hence
\begin{equation*}
\Gamma(a_ia_jt^{n-1};p,q)\Gamma(a_ka_lt^{n-1};p,q)=1.  
\end{equation*}
This implies 
\begin{lem} \label{lem:limJn}
In the limit as $a_2\to a_1^{-1}$, we have 
\begin{equation*}
\begin{split}
&\lim_{a_2\to a_1^{-1}}(1-a_1a_2)J_n(a_1,\ldots,a_6)
\\
&=
\frac{\prod_{i=1}^{n-1}
\Gamma(t^{i};p,q)}{(p;p)_\infty(q;q)_\infty}
\prod_{i=1}^{n}
\prod_{m=3}^{6}
\Gamma(a_1^{\pm1}a_kt^{i-1};p,q)
\prod_{i=1}^{n-1}
\prod_{3\le j<k\le 6}{}
\Gamma(a_ja_kt^{i-1};p,q). 
\end{split}
\end{equation*}
\qed
\end{lem}

We next investigate the behavior of $I_n(a_1,\ldots,a_6)$ as $a_2\to a_1^{-1}$, 
assuming that $|p|$ is sufficiently small so that equality \eqref{eq:I=cJ} holds.  
Here we suppose that $|p|<|q|^7 |t|^{2n-2}$ for convenience.  
As we remarked in Lemma \ref{lem:analcont}, 
in the region \eqref{eq:defUext} 
the integral $I_n(a_1,\ldots,a_6)$ 
is expressed as the integral 
\begin{equation}\label{eq:InInt}
I_n(a_1,\ldots,a_6)=\int_{C^n} \Psi_n(z)\varpi_n(z)
\end{equation}
over a certain $n$-cycle $C^n$, provided that $|t|<r^2$.  
Setting $r=|q|^{\frac{1}{2}}$, we assume further 
\begin{equation*}
1<|a_1|<|q|^{-\frac{1}{2}};\quad |a_m|<1\quad (m=2,\ldots,6).  
\end{equation*}
In this case we can choose the cycle $C$ as 
\begin{equation*}
C=C_0+C_1,\ \ C_1=C_1^{+}-C_1^{-};\quad 
C_0=C_1(0),\quad
C_1^+=C_\epsilon(a_1),\quad C_1^-=C_\epsilon(a_1^{-1}). 
\end{equation*}
Then we analyze the effect of pinching about the cycles $C_1^+$, $C_1^-$ 
as $a_2\to a_1^{-1}$.  

We consider the integral
\begin{equation*}
\frac{1}{2\pi\sqrt{-1}}\int_{C}
\Psi_n(z_1,z_2,\ldots,z_n)\frac{dz_1}{z_1}
\end{equation*}
with respect to $z_1$.  
Since 
\begin{equation*}
\Psi_n(z_1,z_2,\ldots,z_n)=
\Psi_{n-1}(z_2,\ldots,z_n)
\frac{\prod_{m=1}^{6}\Gamma(a_mz_1^{\pm1};p,q)}
{\Gamma(z_1^{\pm2};p,q)}
\prod_{k=2}^{n}
\frac
{\Gamma(tz_1^{\pm1}z_k^{\pm1};p,q)}
{\Gamma(z_1^{\pm1}z_k^{\pm1};p,q)}, 
\end{equation*}
the poles 
$z_1=a_1,a_1^{-1}$ of the integrand 
arise only in the factor
\begin{equation*}
\Gamma(a_1z_1^{\pm1};p,q)
=
\frac
{(pqa_1^{-1}z_1;p,q)_\infty(pqa_1^{-1}z_1^{-1};p,q)_\infty}
{(a_1z_1;p,q)_\infty(a_1z_1^{-1};p,q)_\infty}. 
\end{equation*}
Note that 
\begin{equation*}
\begin{split}
&\mbox{Res}\Big(
\Gamma(a_1z_1^{\pm1};p,q)\frac{dz_1}{z_1};z_1=a_1\Big)
=
\frac
{(pqa_1^{-2};p,q)_\infty}
{(a_1^2;p,q)_\infty
(p;p)_\infty
(q;q)_\infty
}
=
\frac
{\Gamma(a_1^2;p,q)}
{(p;p)_\infty
(q;q)_\infty}
\end{split}
\end{equation*}
and
\begin{equation*}
\mbox{Res}\Big(
\Gamma(a_1z_1^{\pm1};p,q)\frac{dz_1}{z_1};z_1=a_1^{-1}\Big)
=
-\frac
{\Gamma(a_1^2;p,q)}
{(p;p)_\infty(q;q)_\infty}. 
\end{equation*}
Hence we have 
\begin{equation*}
\begin{split}
&\frac{1}{2\pi\sqrt{-1}}\int_{C_1^{\pm}}
\Psi_n(z_1,\ldots,z_n)\dfrac{dz_1}{z_1}
\\
&=
\pm\frac
{\Gamma(a_1^2;p,q)}
{(p;p)_\infty
(q;q)_\infty}
\frac{
\prod_{m=2}^{6}\Gamma(a_ma_1^{\pm1};p,q)}
{\Gamma(a_1^{\pm2};p,q)}
\prod_{k=2}^{n}
\frac
{\Gamma(ta_1^{\pm1}z_k^{\pm1};p,q)}
{\Gamma(a_1^{\pm1}z_k^{\pm1};p,q)}
\Psi_{n-1}(z_2,\ldots,z_n)
\\
&=
\pm\frac
{1}
{(p;p)_\infty
(q;q)_\infty}
\frac{
\prod_{m=2}^{6}\Gamma(a_ma_1^{\pm1};p,q)}
{\Gamma(a_1^{-2};p,q)}
\prod_{k=2}^{n}
\frac
{\Gamma(ta_1^{\pm1}z_k^{\pm1};p,q)}
{\Gamma(a_1^{\pm1}z_k^{\pm1};p,q)}
\Psi_{n-1}(z_2,\ldots,z_n)
\\
&=
\pm\frac
{1}
{(p;p)_\infty
(q;q)_\infty}
\frac{
\prod_{m=2}^{6}\Gamma(a_ma_1^{\pm1};p,q)}
{\Gamma(a_1^{-2};p,q)}
\widehat{\Psi}_{n-1}(z_2,\ldots,z_n), 
\end{split}
\end{equation*}
where 
\begin{equation*}
\begin{split}
&\widehat{\Psi}_{n-1}(z_2,\ldots,z_n)
=
\prod_{k=2}^{n}
\frac
{\Gamma(ta_1^{\pm1}z_k^{\pm1};p,q)}
{\Gamma(a_1^{\pm1}z_k^{\pm1};p,q)}
\Psi_{n-1}(z_2,\ldots,z_n)
\\
&=\prod_{i=2}^{n}
\frac{
\Gamma(ta_1^{\pm1}z_i^{\pm1};p,q)
\prod_{m=2}^{6}\Gamma(a_mz_i^{\pm1};p,q)
}
{\Gamma(a_1^{-1}z_i^{\pm1};p,q)
\Gamma(z_i^{\pm2};p,q)}
\prod_{2\le j<k\le n}
\frac
{\Gamma(tz_j^{\pm1}z_k^{\pm1};p,q)}
{\Gamma(z_j^{\pm1}z_k^{\pm1};p,q)}. 
\end{split}
\end{equation*}
This implies 
\begin{equation}\label{eq:416}
\begin{split}
\frac{1}{2\pi\sqrt{-1}}\int_{C_1}\Psi_n(z)
\dfrac{dz_1}{z_1}
&=
\frac{1}{2\pi\sqrt{-1}}\int_{C_1^+}\Psi_n(z)
\dfrac{dz_1}{z_1}
-
\frac{1}{2\pi\sqrt{-1}}\int_{C_1^-}\Psi_n(z)
\dfrac{dz_1}{z_1}
\\
&=
\frac{
2\prod_{m=2}^{6}\Gamma(a_ma_1^{\pm1};p,q)}
{(p;p)_\infty(q;q)_\infty
\Gamma(a_1^{-2};p,q)}
\widehat{\Psi}_{n-1}(z_2,\ldots,z_n), 
\end{split}
\end{equation}
and hence
\begin{equation*}
\begin{split}
\frac{1}{2\pi\sqrt{-1}}\int_{C}\Psi_n(z)
\dfrac{dz_1}{z_1}
=
\frac{1}{2\pi\sqrt{-1}}\int_{C_0}\Psi_n(z)
\dfrac{dz_1}{z_1}
+
\frac{
2\prod_{m=2}^{6}\Gamma(a_ma_1^{\pm1};p,q)}
{(p;p)_\infty(q;q)_\infty
\Gamma(a_1^{-2};p,q)}
\widehat{\Psi}_{n-1}(z_2,\ldots,z_n).  
\end{split}
\end{equation*}
We remark that the first term is regular at 
$a_1a_2=1$ and has a finite limit as
$a_2\to a_1^{-1}$,
while the second term diverges 
in the order $(1-a_1a_2)^{-1}$ 
because of the factor $\Gamma(a_2a_1;p,q)$. 
Since 
\begin{equation*}
\begin{split}
\frac{2
\Gamma(a_2a_1^{\pm1};p,q)
\prod_{m=3}^{6}
\Gamma(a_ma_1^{\pm1};p,q)}
{(p;p)_\infty(q;q)_\infty
\Gamma(a_1^{-2};p,q)}
&=\frac{1}{1-a_1a_2}
\frac{(pqa_2^{-1}a_1^{-1};p,q)}{(pa_2a_1;p)_\infty
(qa_2a_1;q)_\infty(pqa_2a_1;p,q)_\infty}
\\
&\quad\times
\frac{2
\Gamma(a_2a_1^{-1};p,q)
\prod_{m=3}^{6}
\Gamma(a_ma_1^{\pm1};p,q)}
{(p;p)_\infty(q;q)_\infty
\Gamma(a_1^{-2};p,q)}, 
\end{split}
\end{equation*}
we have 
\begin{equation*}
\begin{split}
\lim_{a_2\to a_1^{-1}}(1-a_1a_2)
\frac{2
\Gamma(a_2a_1^{\pm1};p,q)
\prod_{m=3}^{6}
\Gamma(a_ma_1^{\pm1};p,q)}
{(p;p)_\infty(q;q)_\infty
\Gamma(a_1^{-2};p,q)}
=
\frac{2\,\prod_{m=3}^{6}
\Gamma(a_ma_1^{\pm1};p,q)}
{(p;p)_\infty^2(q;q)_\infty^2}, 
\end{split}
\end{equation*}
where 
$a_6=pq/(a_3a_4a_5t^{2n-2})$
in the right hand side. 
On the other hand, 
\begin{equation*}
\lim_{a_2\to a_1^{-1}}\widehat{\Psi}_{n-1}(z_2,\ldots,z_n)
=
\Psi_{n-1}(ta_1,ta_1^{-1},a_3,\ldots,a_6;z_2,\ldots,z_n).  
\end{equation*}
To summarize, we obtain 
\begin{equation}\label{eq:limintz1}
\begin{split}
&\lim_{a_2\to a_1^{-1}}(1-a_2a_1)
\frac{1}{2\pi\sqrt{-1}}\int_{C}\Psi_n(z_1,\ldots,z_n)
\dfrac{dz_1}{z_1}
\\
&=
\frac{2\,\prod_{m=3}^{6}
\Gamma(a_ma_1^{\pm1};p,q)}
{(p;p)_\infty^2(q;q)_\infty^2}\,
\Psi_{n-1}(ta_1,ta_1^{-1},a_3,\ldots,a_6;z_2,\ldots,z_n). 
\end{split}
\end{equation}

We decompose the multiple integral $I_n(a_1,\ldots,a_6)$ of \eqref{eq:InInt} 
as 
\begin{equation}\label{eq:decomp}
\begin{split}
\int_{C^n}\Psi_n(z)\varpi_ n(z)
&=
\int_{C_1\times C^{n-1}}\Psi_n(z)\varpi_n(z)
+\int_{C_0\times C^{n-1}}\Psi_n(z)\varpi_n(z)
\\[3pt]
&=
\int_{C_1\times C^{n-1}}\Psi_n(z)\varpi_n(z)
+\int_{C_0\times C_1\times C^{n-2}}\Psi_n(z)\varpi_n(z)
\\
&\hspace{40pt}
+\int_{C_0\times C_0\times C^{n-2}}\Psi_n(z)\varpi_n(z)
\\
&=\cdots
\\[5pt]
&=\sum_{i=1}^{n}
\int_{C_0^{i-1}\times C_1\times C^{n-i}}
\Psi_n(z)\varpi_n(z)
+\int_{C_0^{n}}\Psi_n(z)\varpi_n(z). 
\end{split}
\end{equation}
Regarding the integral
\begin{equation*}
\begin{split}
\int_{C_0^{i-1}\times C_1\times C^{n-i}}\Psi_n(z)\varpi_n(z)
&=
\int_{C_0^{i-1}\times C^{n-i}}
\left(
\frac{1}{2\pi\sqrt{-1}}\int_{C_1}
\Psi_n(z)
\frac{dz_i}{z_i}
\right)
\varpi_{n-1}(z_{\,\widehat{i}\,}), 
\\[2pt]
\varpi_{n-1}(z_{\,\widehat{i}\,})
&=
\frac{1}{(2\pi\sqrt{-1})^{n-1}}
\prod_{1\le j\le n\atop j\ne i}{}\frac{dz_j}{z_j}, 
\end{split}
\end{equation*}
where $z_{\,\widehat{i}}=
(z_1,\ldots,z_{i-1},z_{i+1},\ldots,z_n)$, by \eqref{eq:416} we have  
\begin{equation*}
\begin{split}
\frac{1}{2\pi\sqrt{-1}}\int_{C_1}
\Psi_n(z)
\frac{dz_i}{z_i}
=
\frac{
2\,\Gamma(a_2a_1^{\pm1};p,q)\,
\prod_{m=3}^{6}\Gamma(a_ma_1^{\pm1};p,q)}
{(p;p)_\infty(q;q)_\infty
\Gamma(a_1^{-2};p,q)}
\widehat{\Psi}_{n-1}(z_{\,\widehat{i}\,}).  
\end{split}
\end{equation*}
Since 
$\widehat{\Psi}_{n-1}(z_{\,\widehat{i}\,})
=\widehat{\Psi}_{n-1}(z_1,\ldots,z_{i-1},z_{i+1},\ldots,z_n)
$ 
is regular at 
$z_j=a_1,a_1^{-1}$ $(j\ne i)$, 
one can replace the $(n-1)$-cycle 
$C_0^{i-1}\times C^{n-i}$ 
by 
$C_0^{n-1}$ as 
\begin{equation*}
\begin{split}
\int_{C_0^{i-1}\times C_1\times C^{n-i}}
\Psi_n(z)\varpi_n(z)
=
\frac{
2\,\Gamma(a_2a_1^{\pm1};p,q)\,
\prod_{m=3}^{6}\Gamma(a_ma_1^{\pm1};p,q)}
{(p;p)_\infty(q;q)_\infty
\Gamma(a_1^{-2};p,q)}
\int_{C_0^{n-1}}
\widehat{\Psi}_{n-1}(z_{\,\widehat{i}\,})
\varpi(z_{\,\widehat{i}\,}).  
\end{split}
\end{equation*}
Hence \eqref{eq:decomp} implies 
\begin{equation}\label{eq:426}
\begin{split}
&
\int_{C^n}\Psi_n(z)\varpi_n(z)
\\
&\quad
=\frac{
2n\Gamma(a_2a_1^{\pm1};p,q)
\prod_{m=3}^{6}\Gamma(a_ma_1^{\pm1};p,q)}
{(p;p)_\infty(q;q)_\infty
\Gamma(a_1^{-2};p,q)}
\int_{C_0^{n-1}}
\widehat{\Psi}_{n-1}(z_{\,\widehat{1}})
\varpi_{n-1}(z_{\,\widehat{1}})
+\int_{C_0^n}\Psi_n(z)\varpi_n(z).  
\end{split}
\end{equation}
Note that the second term of the right-hand side has 
a finite limit as $a_2\to a_1^{-1}$.  
Multiplying \eqref{eq:426} by $1-a_1a_2$, 
we compute the limit 
$a_2\to a_1^{-1}$ by \eqref{eq:limintz1}
as 
\begin{equation*}
\begin{split}
&\lim_{a_2\to a_1^{-1}}(1-a_1a_2)
I_n(a_1,\ldots,a_6)
=
\lim_{a_2\to a_1^{-1}}(1-a_1a_2)
\int_{C^n}\Psi_n(z)\varpi_n(z)
\\
&=
\frac{
2n\prod_{m=3}^{6}\Gamma(a_ma_1^{\pm1};p,q)}
{(p;p)_\infty^2(q;q)_\infty^2}
\int_{C_0^{n-1}}
\Psi_{n-1}(ta_1,ta_1^{-1},a_3,\ldots,a_6;z_{\,\widehat{1}})
\varpi_{n-1}(z_{\,\widehat{1}}) 
\\
&=
\frac{
2n\prod_{m=3}^{6}\Gamma(a_ma_1^{\pm1};p,q)}
{(p;p)_\infty^2(q;q)_\infty^2}
I_{n-1}(ta_1,ta_1^{-1},a_3,\ldots,a_6).
\end{split}
\end{equation*}
This means that 
\begin{equation*}
\begin{split}
&c_n\lim_{a_2\to a_1^{-1}}(1-a_1a_2)J_n(a_1,\ldots,a_6)
\\
&=
\frac{
2n\prod_{m=3}^{6}\Gamma(a_ma_1^{\pm1};p,q)}
{(p;p)_\infty^2(q;q)_\infty^2}\,
c_{n-1}\,J_{n-1}(ta_1,ta_1^{-1},a_3,\ldots,a_6). 
\end{split}
\end{equation*}
The left-hand side is already given by Lemma \ref{lem:limJn}, 
while 
the right-hand side is computed as 
\begin{equation*}
\begin{split}
\frac{2n\,c_{n-1}}{(p;p)_\infty^2(q;q)_\infty^2}
\prod_{i=2}^{n}
\Gamma(t^{i};p,q)
\prod_{i=1}^{n}
\prod_{k=3}^{6}
\Gamma(a_1^{\pm1}a_kt^{i-1};p,q)
\prod_{i=1}^{n-1}
\prod_{3\le j<k\le 6}
\Gamma(a_ja_kt^{i-1};p,q)
\end{split}
\end{equation*}
by definition \eqref{eq:defJ}. 
Comparing these two expressions we obtain the 
recurrence formula 
\begin{equation*}
c_n=c_{n-1}
\frac{2n}{(p;p)_\infty(q;q)_\infty}
\frac{\Gamma(t^{n};p,q)}{\Gamma(t;p,q)}
\qquad(n=1,2,\ldots)
\end{equation*}
for the constants  $c_n$. 
Starting from $c_0=1$, we have 
\begin{equation*}
c_n=
\frac{2^nn!}{(p;p)_\infty^n(q;q)_\infty^n}
\prod_{i=1}^{n}\frac{\Gamma(t^i;p,q)}{\Gamma(t;p,q)}
\qquad(n=0,1,2,\ldots).  
\end{equation*}
This completes the evaluation of the $BC_n$ elliptic Selberg integral 
\begin{equation*}
\begin{split}
I_n(a_1,\ldots,a_6)&=c_n J_n(a_1,\ldots,a_6)
\\
&
=\frac{2^nn!}{(p;p)_\infty^n(q;q)_\infty^n}
\prod_{i=1}^{n}\frac{\Gamma(t^i;p,q)}{\Gamma(t;p,q)}\,
\prod_{i=1}^{n}\,
\prod_{1\le j<k\le 6}\Gamma(a_ja_kt^{i-1};p,q).  
\end{split}
\end{equation*}

\section*{Acknowledgements} 
This work is supported by JSPS Kakenhi Grants (C)25400118 and 
(B)15H03626. 


{\small
\begin{thebibliography}{99}

\bibitem{Ao}
K.~Aomoto,
Jacobi polynomials associated with Selberg integrals, SIAM J. Math. Anal. \textbf{18} (1987), 545--549. 

\bibitem{An} 
G.~W.~Anderson, 
A short proof of Selberg's generalized beta formula, 
Forum. Math. \textbf{3} (1991), 415--417.
%
\bibitem{CG2007}
H.~Coskun and R.~A.~Gustafson, Well-poised Macdonald functions $W_\lambda$ and Jackson coefficients $\omega_\lambda$ on $BC_n$; 
in {\em Jack, Hall--Littlewood and Macdonald polynomials}, pp.127--155, Contemp. Math., {\bf 417}, Amer. Math. Soc., Providence, RI, 2006. 
%
\bibitem{Di}
A.~L.~Dixon, Generalizations of Legendre's formula $K E' - (K - E) K' =
{1 \over 2} \pi$, Proc. London Math. Soc. \textbf{3} (1905), 206--224.

\bibitem{F}
P.~J.~Forrester, Log-Gases and Random Matrices, Princeton University Press, Princeton, 2010.

\bibitem{Gu}
R.~A.~Gustafson, A generalization of Selberg's beta integral, 
Bull. Amer. Math. Soc. (N.S.) \textbf{22} (1990), 97--105.

\bibitem{I}
M.~Ito, Three-term relations between interpolation polynomials for a $BC_n$-type basic hypergeometric series, Adv. Math. \textbf{226} (2011), 4096--4130.

\bibitem{IF1}
M.~Ito and P.~J.~Forrester, 
The $q$-Dixon--Anderson integral and multi-dimensional ${}_1\psi_1$ summations, J. Math. Anal. Appl. \textbf{423} (2015), 1704--1737.

\bibitem{IF2} 
M.~Ito and P.~J.~Forrester, 
A bilateral extension of the $q$-Selberg integral, 
%
Trans. Amer. Math. Soc., to appear.

\bibitem{IN}
M.~Ito and M.~Noumi,
Derivation of a $BC_n$ elliptic summation formula
via the fundamental invariants, 
%
Constr.\,Approx., to appear.
%

\bibitem{IN15}
M.~Ito, M.~Noumi, A generalization of the Sears--Slater transformation 
and elliptic Lagrange interpolation of type $BC_n$, 
Adv. Math. \textbf{299} (2016), 361--380.
%

\bibitem{NR}
B.~Nassrallah and M.~Rahman, 
Projection formulas, a reproducing kernel and a generating function for $q$-Wilson polynomials, 
SIAM J. Math. Anal. \textbf{16} (1985), 186--197.
%
\bibitem{Ra2006}
E.~M.~Rains, $BC_n$-symmetric Abelian functions. Duke Math. J. {\bf 135} (2006), 99--180.

\bibitem{R}
E.~M.~Rains, Transformations of elliptic hypergeometric integrals, Ann. of Math. (2) \textbf{171} (2010), 169--243.

\bibitem{Se}
A.~Selberg, Remarks on a multiple integral, Norsk Mat. Tidsskr. \textbf{26} (1944). 71--78. 

\bibitem{S}
V.~P.~Spiridonov, Short proofs of the elliptic beta integrals, Ramanujan J. \textbf{13} (2007), 265--283. 

\bibitem{vDS}
J.~F.~van Diejen and V.~P.~Spiridonov, Elliptic Selberg integrals, Internat. Math. Res. Notices \textbf{20} (2001), 1083--1110.

\end{thebibliography}
}
\end{document}